# A Math Circle in Third Grade

## A math circle in a third-grade primary school

**Evgeny Lakshtanov** (University of Aveiro, lakshtanov@ua.pt), **Marek Tusiewicz** (marek@tusiewicz.eu)

*Market-game spreadsheet template:* https://github.com/lakshtanov/kolo-matematyczne

---

## 1. Introduction

Over the course of a school year, math circle sessions were held in a third-grade class at a Polish primary school. The content was developed by Evgeny Lakshtanov — a mathematician who speaks Polish only at a basic level. The sessions were led by Marek Tusiewicz. Marek is not a professional mathematician, but he has a deep interest in details and a rare ability to communicate with children — he could sense when the class was losing interest, when the explanation was too fast, when a pause was needed. Evgeny's language barrier shaped the entire format: minimal explanations, maximum action.

The school was an ordinary primary school with one teacher per class. Both authors' daughters attended the same class — Marta (Evgeny's daughter) and Iga (Marek's daughter). The authors would bring their children to school in the morning, step into the classroom for 45 minutes (one lesson, until the bell), and then go to work. The **entire class** participated — about 25 children aged 9–10 — not just those who showed particular interest in mathematics.

The evening before each session, the authors would meet: they discussed the lesson plan in detail, solved the problems together, and prepared backup activities in case something became boring. After each lesson, they discussed the didactic aspects — what worked, what didn't, where the children lost the thread — and jointly decided whether the plan needed adjustment for the next session.

## 2. Methodological Foundations

Henri Poincaré in *Science and Method* (1908) described the mechanism of mathematical discovery: conscious work on a problem gives way to a period of rest, during which the subconscious continues to combine ideas. Then — a sudden illumination, when a ready solution "surfaces." But the subconscious can only work with the material that consciousness has provided: there must be a preliminary phase of saturation — intensive, even if fruitless, work on the problem.

This principle applies to teaching children. It is important for children to encounter complex mathematical ideas as early as possible — not in a form simplified to meaninglessness, but in a form that preserves intellectual value. Familiarity should be reinforced through repeated exercises: playing the same game with variations saturates the subconscious with material. A child cannot formulate a strategy, but after twenty games starts to "feel" the right move — the subconscious has processed the experience and produced intuition.

That is why in our circle we do not rush to explain. First — many games. Then, when children themselves begin to notice patterns and want to discuss them — a brief explanation. Not before. Premature explanation kills the process of independent discovery and deprives the subconscious of its work.

This approach resonates with the tradition of math circles: Zvonkin's home circle for preschoolers [1] showed that children aged 4–7 can work with combinatorics and graph theory through play. The question of when and how a complex theory can be simplified so that the simplification remains substantive is explored in [2]. The general principle of early exposure to deep ideas through problems is the foundation of the Kolmogorov–Arnold–MCCME circle tradition [3].

## A heterogeneous class and the zone of maximal development

Our circle worked with the entire class — a heterogeneous group. Some children calculate quickly, others slowly; some see patterns immediately, others need ten games. The question arises: how to choose the right difficulty level?

Imagine a child's rate of development as a function of task difficulty. This function is not monotonic — it has a maximum:

- **Task too easy** (everything is clear, can do it alone) — rate of development is low. The child repeats what they already know. No challenge — no growth.
- **Task in the right zone** (roughly understand the idea, but cannot produce it independently yet) — rate of development is maximal. The child stretches, tries, makes mistakes, gropes — this is precisely where understanding forms.
- **Task too hard** (completely incomprehensible, nothing to grasp) — rate of development is again low or zero. The child disconnects.

This is essentially a gradient: the rate of development is the derivative of understanding with respect to time, and it is maximal in the intermediate zone. Vygotsky called this the **zone of proximal development** [5] — the region between what a child can do alone and what they can do with help. Zvonkin [1] discusses this effect in detail using examples from his circle.

The problem is that in a heterogeneous class, this zone is different for each child. We cannot select the ideal difficulty for everyone simultaneously. But we can try to keep the **majority** of children near their maximum — so that the average rate of development of the class is as high as possible.

The competitive format helps: strong children don't get bored (they compete or help the team), weak children don't disconnect (excitement maintains attention, the team helps). Repeated play of the same game with variations allows everyone to progress at their own pace — some master the strategy on the fifth game, some on the twentieth.

**Disclaimer.** When below we say that "children understood" a strategy or "saw" an isomorphism, this should be understood through the lens of the approach described above. Understanding is heterogeneous: some children indeed formulated the strategy in words, others only intuitively felt the right move, still others were only observing — but all were in the process, each at their own level.

Our circle implements these principles for third grade: NIM is combinatorial game theory, the market is stochastic optimization, weights are the ternary system, guess the number is information theory. Children don't know these words, but they work with the ideas themselves.

## 3. Goals

The sessions had two equal goals:

**Goal 1: introduce mathematical ideas.** Strategy, optimization, the concept of equality, binary search, backward analysis — through games, without formal definitions and without lectures.

**Goal 2: make them calculate.** Arithmetic not as an end in itself, but as a tool of the game. To win at "Market," you need to quickly multiply and subtract. To win at NIM — count remainders. The motivation comes from the desire to win, not from external authority.

## 4. Format: Three Teams and Competition

With the teacher's help, the class is divided into **three balanced teams** (by rows). Teams play in pairs — two compete, the third observes. Two students from each team come to the board. The authors monitored rotation: **all children participated**, not just the active or strong students.

Competition is the main mechanism for maintaining attention. It solves several problems at once:

- **Focus.** Children don't get distracted because they're watching the score.
- **Motivation to calculate.** Don't want to let the team down — calculate correctly and quickly.
- **Social skills.** Must listen to others on the team, negotiate, make collective decisions.

Explanations are kept to a minimum — brief, only when children themselves want to listen. The main experience comes through repeated play, gradually discovering strategies.

---

## 5. Game 1: Expressions with a Die

*Game with expressions on the board. Frequently at the beginning of the year — the main starting game.*

### 5.1. Rules

An arithmetic expression with blank squares is written on the board — identical for both players. A die is rolled. Each player decides which square to fill with the rolled number. Whoever gets the larger value of the expression at the end wins.

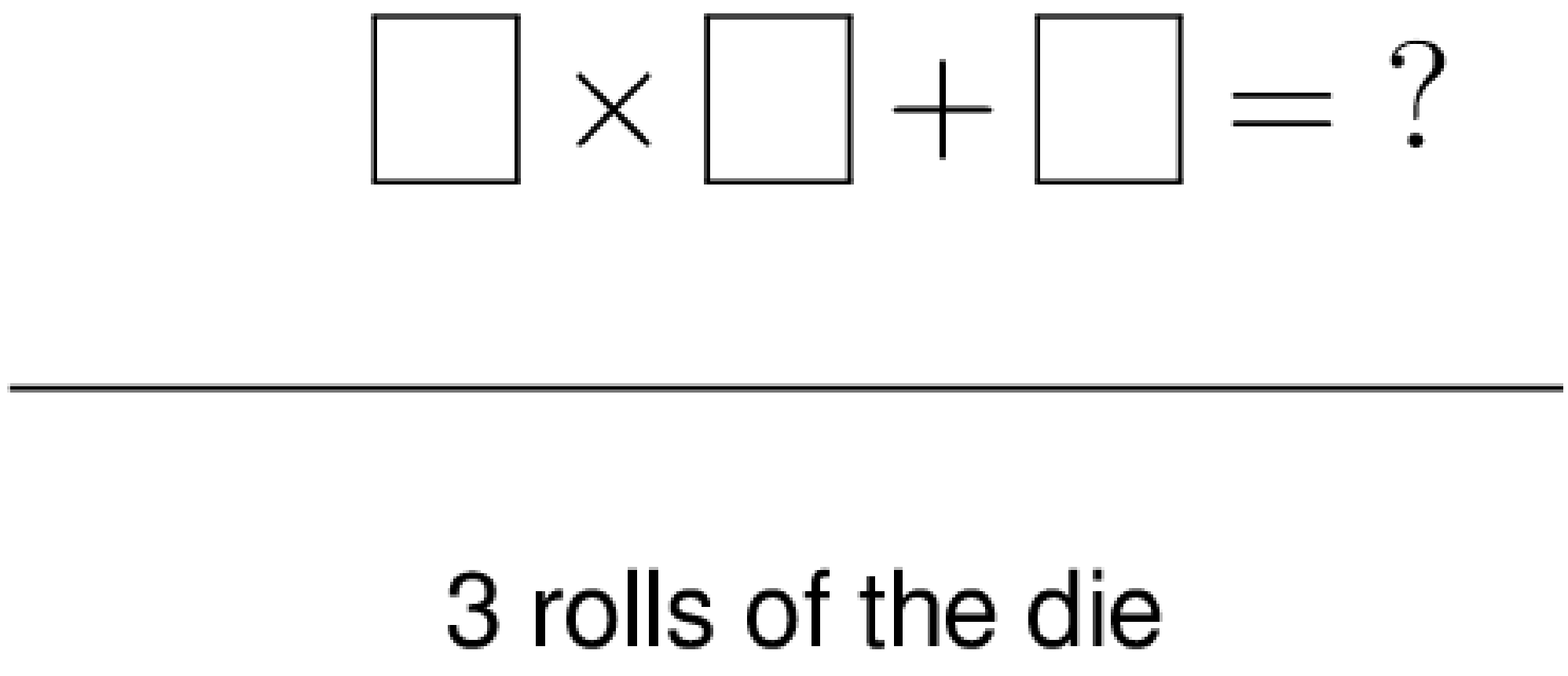

*Fig. 1. Expression*

### 5.2. Example: expression with strategy

On the board: ▢ × ▢ + ▢. Three rolls, three squares. If the goal is maximum, large numbers are better placed in the multipliers, not in the addend.

Roll 1: got **6**. Where to put it? In a multiplier — potential contribution up to 6 × 6 = 36. In the addend — just +6. Strategy: **large numbers go in multipliers**.

Suppose the rolls are: 6, 2, 4. Optimal placement: 6 × 4 + 2 = 26. Poor: 2 × 4 + 6 = 14.

### 5.3. Expression variants

Throughout the year we used different expressions, each with its own strategy:

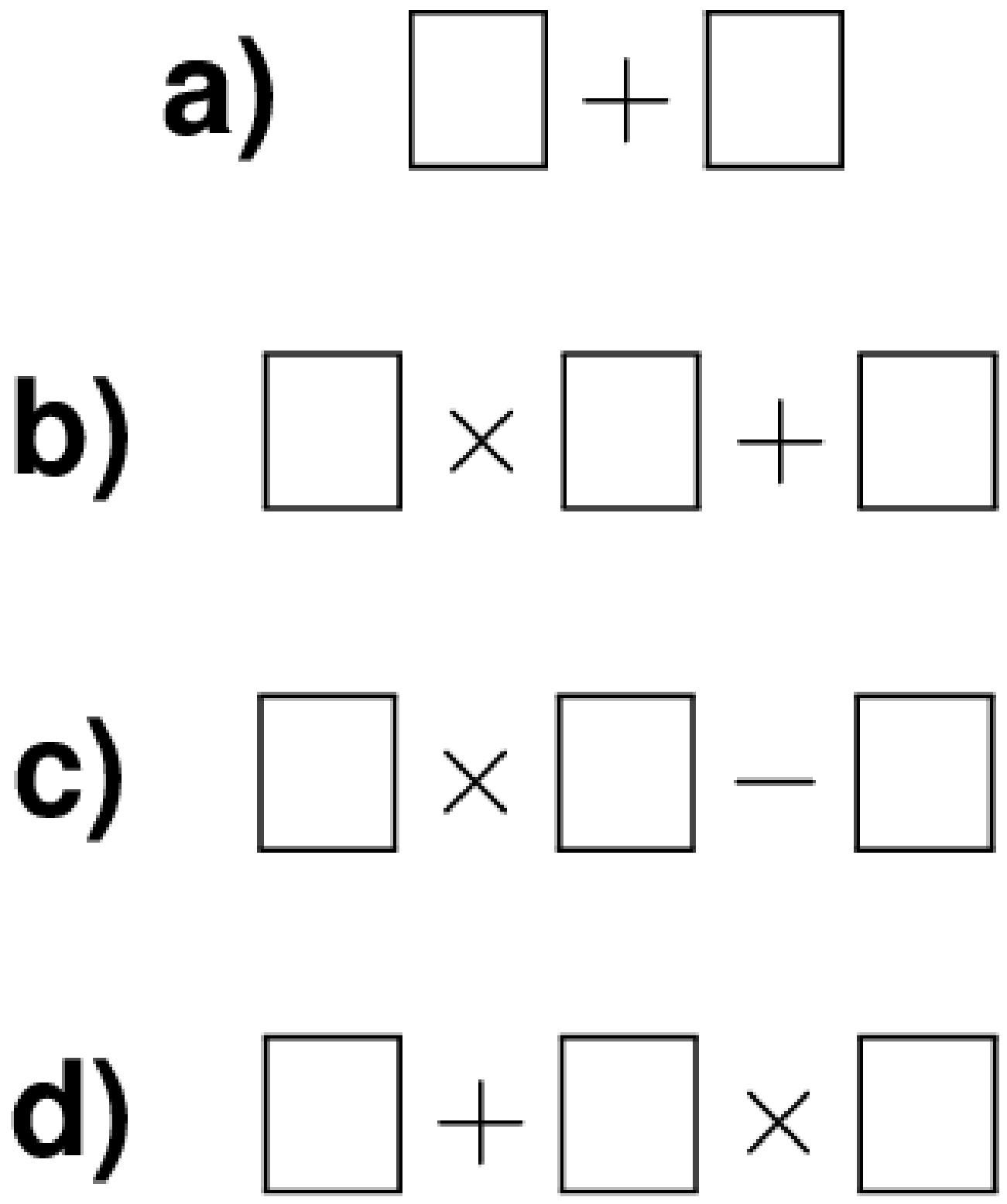

*Fig. 2. Expression variants*

**a)** ▢ + ▢ — simple sum, no strategy (commutativity). Warm-up.

**b)** ▢ × ▢ + ▢ — main expression. Large numbers are better in the multiplier (die average ≈ 3.5).

**c)** ▢ × ▢ − ▢ — subtraction changes the strategy: small numbers are better subtracted.

**d)** ▢ + ▢ × ▢ — order of operations! Multiplication goes first. Children who don't know the priority lose → they learn.

**e)** (▢ + ▢) × ▢ — brackets change the strategy: now the sum is multiplied as a whole, and a large number may be better inside the brackets.

**f)** ▢ × ▢ + ▢ × ▢ — four rolls, two products. More complex evaluation.

When the expression is slightly modified — the strategy changes too. This teaches **flexibility of thinking**, not recipe memorization.

### 5.4. Rule modifications

**Dice:**

• Standard die (faces 1–6)
• Sum of two dice (range 2–12, distribution not uniform — this was also discussed)
• 20-sided die, d20 (range 1–20, greater spread — harder strategy)

**Rules:**

• Maximum → minimum (whoever gets **less** wins). We sometimes changed the rules mid-series: played several rounds for maximum, then announced "now minimum!" — the strategy flips, and children had to adapt.
• Instead of one expression — **two to choose from**. The player circles which expression to fill in.
• **Hit the target**: the goal is not maximum, but to get closest to a given number (e.g., 50). The strategy changes completely: now large numbers aren't always good. Children learn to **estimate**, not just maximize.

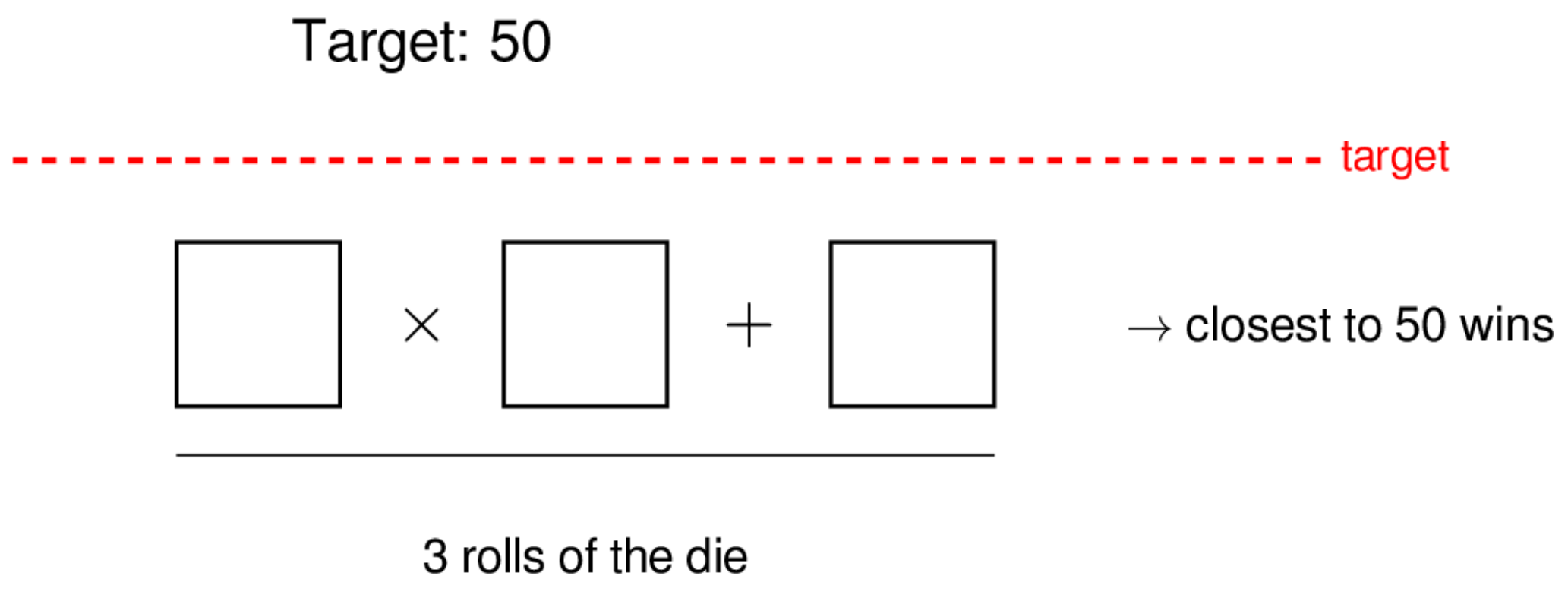

*Fig. 3. Hit the target*

## 6. Game 2: Apple Market

*Apple market. Regular throughout the year, one of the main and most emotional games.*

## 6.1. Rules

Each team starts with **100 złoty** and **0 apples**. The initial apple price is 10 złoty. Each turn the price randomly changes by ±1 or ±2. The team can buy or sell any number of apples at the current price. After a set number of turns, the team with the most **money** (not apples!) wins.

## 6.2. How it looked in class

On the projector screen — Marek's Excel spreadsheet. Everything is visible: current price, how many apples and how much money each team has, history of recent moves. Marek announces the new price aloud and simultaneously it appears on screen.

**Projector screen – Marek's table (turns 4–12)**

| | | 100 Actions | | | Mathematics | | | Lizards | | |
|---|---|---|---|---|---|---|---|---|---|---|
| Turn | Price | Buy/Sell | Apples | Money | Buy/Sell | Apples | Money | Buy/Sell | Apples | Money |
| 4 | 8 zł | 2 | 6 | 47 zł | 1 | 7 | 39 zł | 1 | 4 | 63 zł |
| 5 | 7 zł | 1 | 7 | 40 zł | 1 | 8 | 32 zł | −1 | 3 | 70 zł |
| 6 | 9 zł | −1 | 6 | 49 zł | −3 | 5 | 59 zł | 2 | 5 | 52 zł |
| 7 | 11 zł | −2 | 4 | 71 zł | −2 | 3 | 81 zł | −1 | 4 | 63 zł |
| 8 | 10 zł | 0 | 4 | 71 zł | 0 | 3 | 81 zł | 0 | 4 | 63 zł |
| 9 | 9 zł | 1 | 5 | 62 zł | 1 | 4 | 72 zł | 0 | 4 | 63 zł |
| 10 | 11 zł | −1 | 4 | 73 zł | −2 | 2 | 94 zł | −3 | 1 | 96 zł |
| 11 | 13 zł | −2 | 2 | 99 zł | −2 | 0 | 120 zł | −1 | 0 | 109 zł |
| 12 | 12 zł | −1 | 1 | 111 zł | 0 | 0 | 120 zł | 0 | 0 | 109 zł |

*Fig. 4. Projector screen: table with price, trades, and each team's state*

Children saw exactly this table — not graphs. The current price, who bought/sold what, how many apples and how much money each team has. Graphs (price over time, portfolio values) were sometimes shown at the end of a game for discussion — but during play, decisions were made from the table.

Teams have **10 seconds** to decide. Within the row — heated whispered discussion: "buy!", "no, wait, the price will drop!", "we don't have much money left!" But only the **broker** — the one holding the stick — announces the decision. They stand and say: "Buying three" or "Selling two" or "Pass."

Marek enters the decision in Excel. Everyone sees the numbers change: one team has less money and more apples. Another — the opposite. The third waits.

Then the stick passes to the next child in the row. New broker for the next turn.

## 6.3. Emotions and discipline

The Market is the most emotional game. Price drops — the team that just bought expensive groans. Price rises — those who sold early regret it. Arguments within teams were heated: one wants to buy, another to wait, a third to sell everything.

Children learned several things simultaneously:

- **Managing emotions:** price dropped — don't panic, think.
- **Listening to others:** you're not alone, the decision is collective.
- **Taking responsibility:** the broker decides for everyone, and if the decision was bad — the team doesn't get angry, but discusses what to do next.
- **Living with uncertainty:** nobody knows the future price, but one can reason.

The broker stick appeared out of necessity: without it, everyone shouted simultaneously, and Marek couldn't hear the decision. With the stick — discipline: one speaks, the rest listen.

## 6.4. When the broker sabotages

A separate story: some children simply liked selling everything. The broker would stand and announce "selling everything!" — regardless of what the team whispered. The team would lose the position accumulated over several turns.

This caused real conflict within teams. And then the authors **suggested children write trading rules** — restrictions for the broker. In essence, they invented market regulation: the broker cannot sell more than N apples per turn, or must consult the team before a large trade. Rules were different for different teams.

This is one of the most valuable pedagogical moments: children faced a problem (abuse of power), discussed it, and found an institutional solution — not punishment, but a rule.

All three teams wrote their rules. Remarkably, children formulated both **strategy** (buy below 10, sell above 10, sell everything at the end) and **social norms** (cooperate, don't argue, suggest but don't command). They don't separate mathematics and ethics — for them these are parts of the same game.

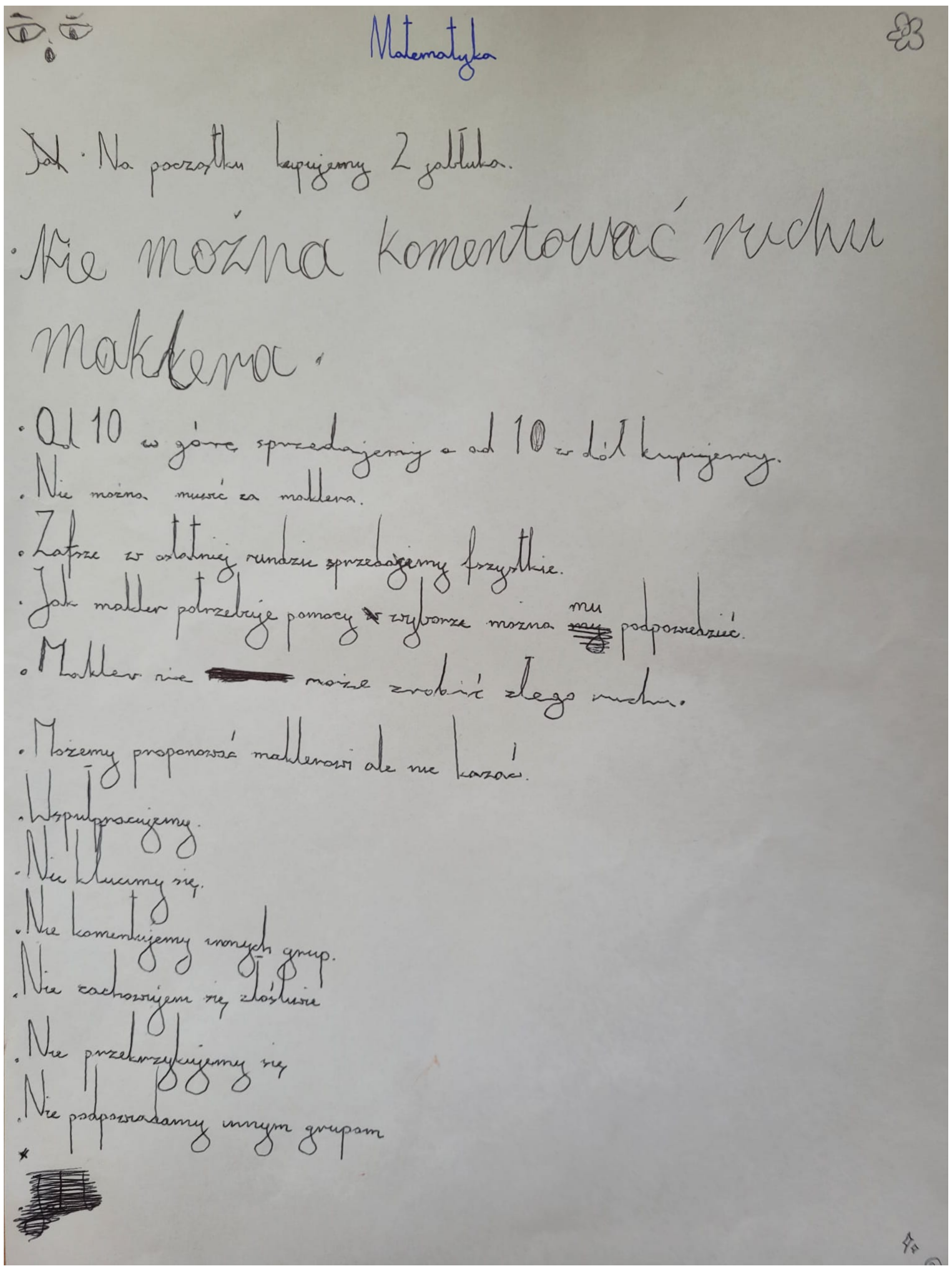
Matematyka

- Na pozatku kupujemy 2 jabłka.
- Nie można komentować ruchu maklera.
- Od 10 w górę sprzedajemy a od 10 w dół kupujemy.
- Nie można musić za maklera.
- Zafsze w ostatniej rundzie sprzedajemy fszytkie.
- Jak makler potrzebuje pomocy w wyborze można mu podpowiedzieć.
- Makler nie [redacted] może zrobić złego ruchu.
- Możemy proponować maklerowi ale nie kazać.
- Wspulpracujemy.
- Nie kłucimy się.
- Nie komentujemy innych grup.
- Nie zachowujem się złośliwie
- Nie przekrzykujemy się
- Nie podpowadamy innym grupam

*Fig. 5. Team 1 rules*

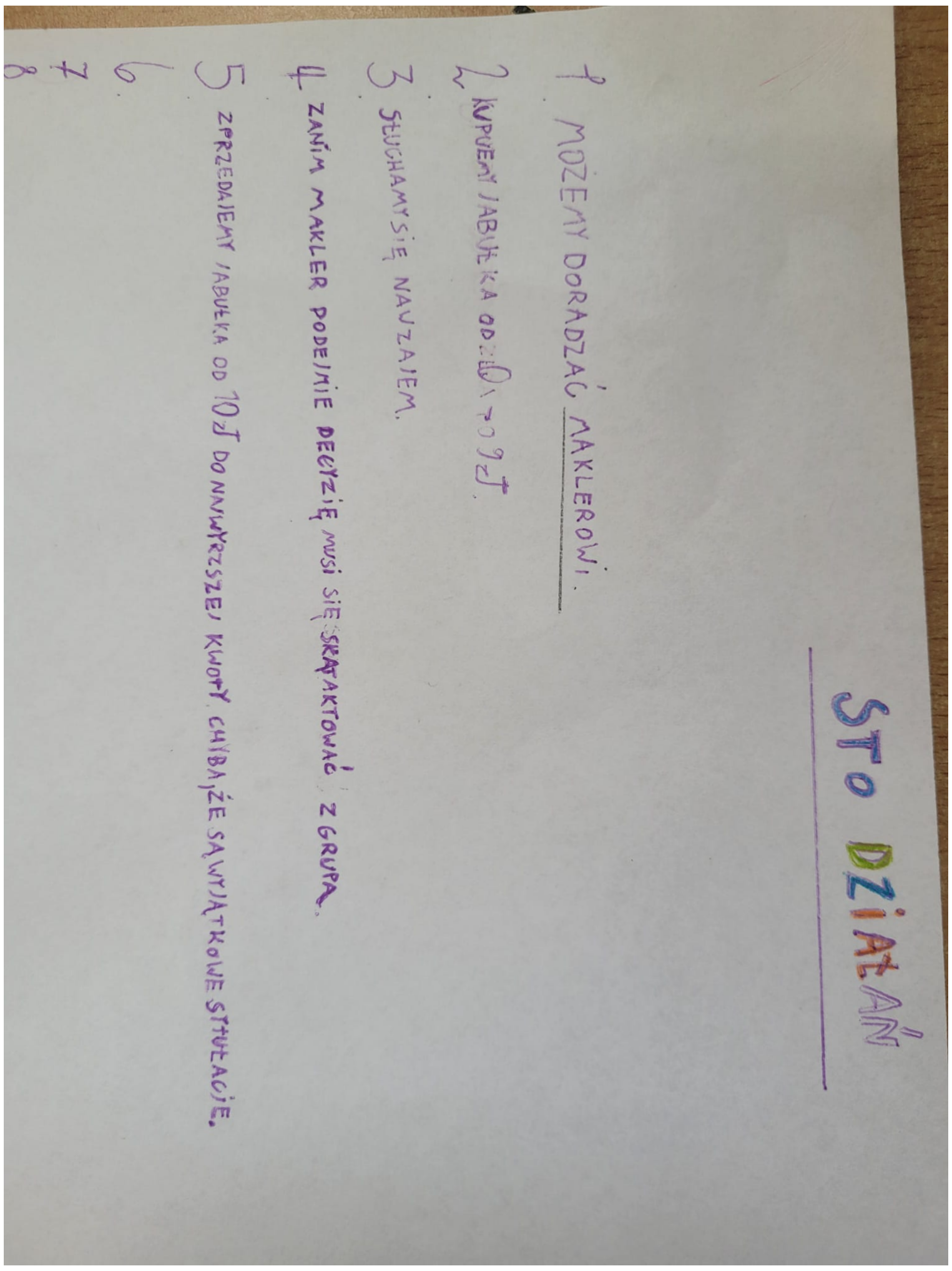

*Fig. 6. Team 2 rules — "Sto działań"*

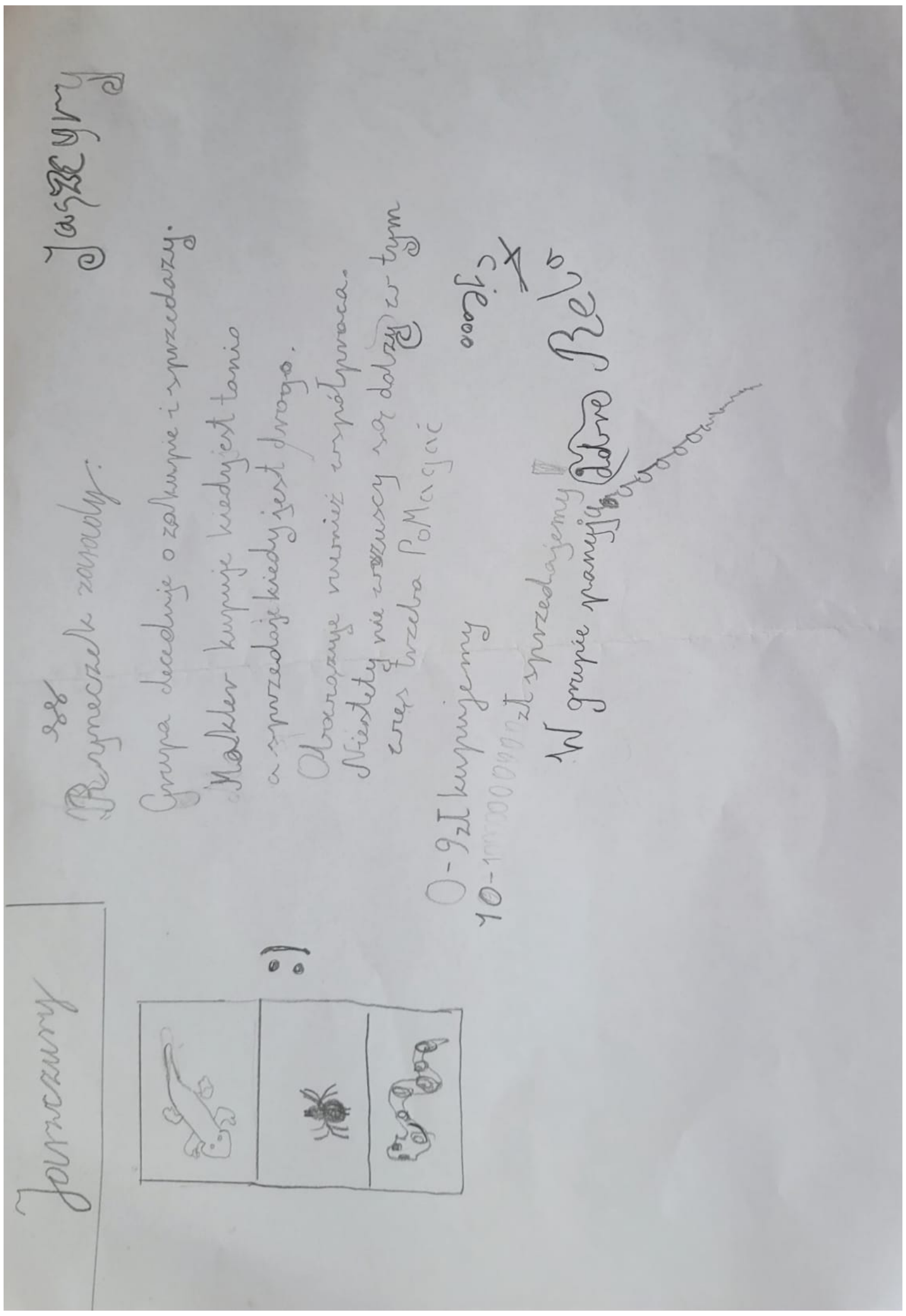

*Fig. 7. Team 3 rules — "Jaszczury" (Lizards)*

## 6.5. Sample game

### Projector screen – Marek's table (turns 4–12)

| | | 100 Actions | | | Mathematics | | | Lizards | | |
|---|---|---|---|---|---|---|---|---|---|---|
| Turn | Price | Buy/Sell | Apples | Money | Buy/Sell | Apples | Money | Buy/Sell | Apples | Money |
| 4 | 8 zł | 2 | 6 | 47 zł | 1 | 7 | 39 zł | 1 | 4 | 63 zł |
| 5 | 7 zł | 1 | 7 | 40 zł | 1 | 8 | 32 zł | –1 | 3 | 70 zł |
| 6 | 9 zł | –1 | 6 | 49 zł | –3 | 5 | 59 zł | 2 | 5 | 52 zł |
| 7 | 11 zł | –2 | 4 | 71 zł | –2 | 3 | 81 zł | –1 | 4 | 63 zł |
| 8 | 10 zł | 0 | 4 | 71 zł | 0 | 3 | 81 zł | 0 | 4 | 63 zł |
| 9 | 9 zł | 1 | 5 | 62 zł | 1 | 4 | 72 zł | 0 | 4 | 63 zł |
| 10 | 11 zł | –1 | 4 | 73 zł | –2 | 2 | 94 zł | –3 | 1 | 96 zł |
| 11 | 13 zł | –2 | 2 | 99 zł | –2 | 0 | 120 zł | –1 | 0 | 109 zł |
| 12 | 12 zł | –1 | 1 | 111 zł | 0 | 0 | 120 zł | 0 | 0 | 109 zł |

*Fig. 8. Sample game*

| Turn | Price | Action | Apples | Money |
|---|---|---|---|---|
| Start | 10 | — | 0 | 100 |
| 1 | 12 | Buy 4 | 4 | 100 − 4×12 = 52 |
| 2 | 11 | Buy 3 | 7 | 52 − 3×11 = 19 |
| 3 | 13 | Sell 2 | 5 | 19 + 2×13 = 45 |
| 4 | 12 | Hold | 5 | 45 |
| 5 | 14 | Sell 5 | 0 | 45 + 5×14 = 115 |

Total: started with 100 złoty, ended with 115. Profit: **+15 złoty**.

## 6.6. Arithmetic as a tool

Each turn requires calculations — and children do them not because they were asked, but because **without calculation, no decision can be made**:

- "We have 52 złoty, an apple costs 13. How many can we buy?" → 52 ÷ 13 = 4.
- "Buying 3 at 11. How much money left?" → 52 − 3 × 11 = 52 − 33 = 19.
- "Selling 5 at 14. How much do we get?" → 5 × 14 = 70.

An arithmetic error is not a bad mark in a notebook, but a real loss for the team. The motivation to calculate correctly becomes intrinsic.

## 6.7. Worksheets

The Market showed that children calculate slowly — multiplication and subtraction of two-digit numbers was difficult. This dissatisfaction led to a separate block of arithmetic lessons (see Chapter 9), during which we discovered a problem with the equals sign and arrived at the idea of the scales lesson.

In a separate session, children received tables with missing numbers — training market arithmetic in a calm setting:

| Turn | Price | Buy(+)/Sell(−) | Apples | Money |
|---|---|---|---|---|
| Start | 10 | 0 | 0 | 100 |
| 1 | 12 | +4 | ? | ? |
| 2 | 11 | +3 | 7 | 19 |
| 3 | 13 | ? | 5 | ? |
| 4 | 12 | −1 | ? | 45 |
| 5 | 14 | ? | 0 | ? |

Solution: after turn 1, apples = 0 + 4 = 4, money = 100 − 4 × 12 = 52. Turn 3: had 7 apples, now 5 → sold 2, money = 19 + 2 × 13 = 45.

## 6.8. Lessons with predetermined scenarios

In one lesson, the price was not random but predetermined — to demonstrate the idea of strategy:

**"Growth" scenario:** 10, 11, 12, 13, 14, 15, 16, 17, 18, 19, 20. Optimal strategy: buy as much as possible on the first turn (at 10), sell on the last (at 20).

**"Decline" scenario:** 10, 9, 8, 7, 6, 5, 4, 3, 2, 1. Optimal strategy: buy nothing.

**"Saw" scenario:** 10, 14, 10, 14, 10, 14, ... Optimal strategy: buy at 10, sell at 14. Profit of 4 złoty per apple every two turns.

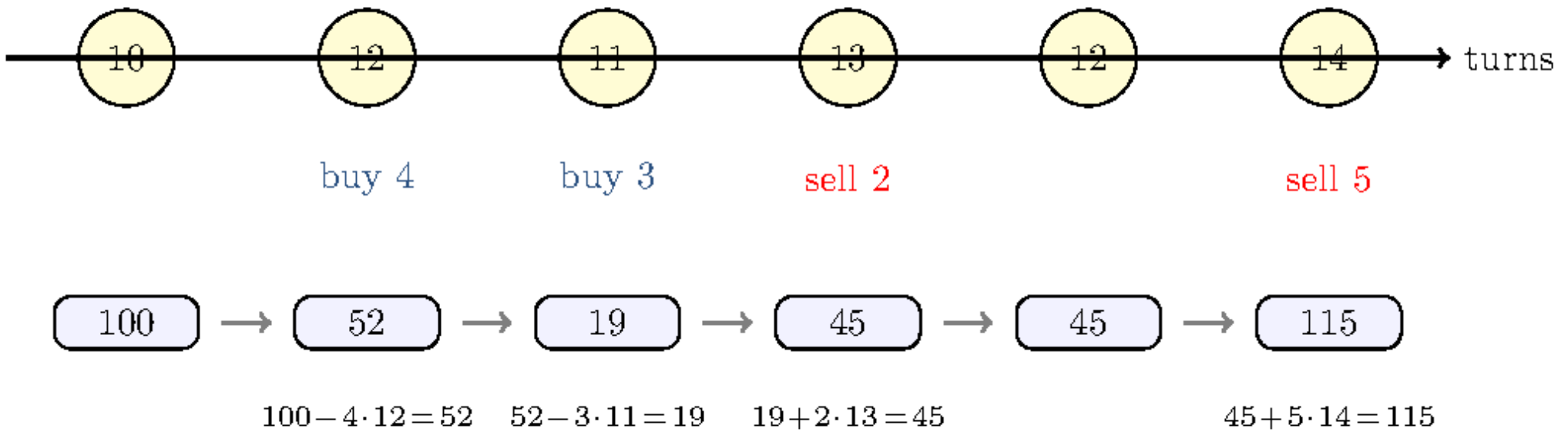

*Fig. 9. Price scenarios*

It was **not obvious** to children what to do, even when the trend was visible! But through repeated games they gradually figured it out — and competition accelerated this process, because the losing team could see that their opponents were doing something right.

---

## 7. Combinatorial Games: NIM and Chodniczek

*A family of games in which children gradually discovered the "backward analysis" strategy and — most importantly — saw that different games can be the same game.*

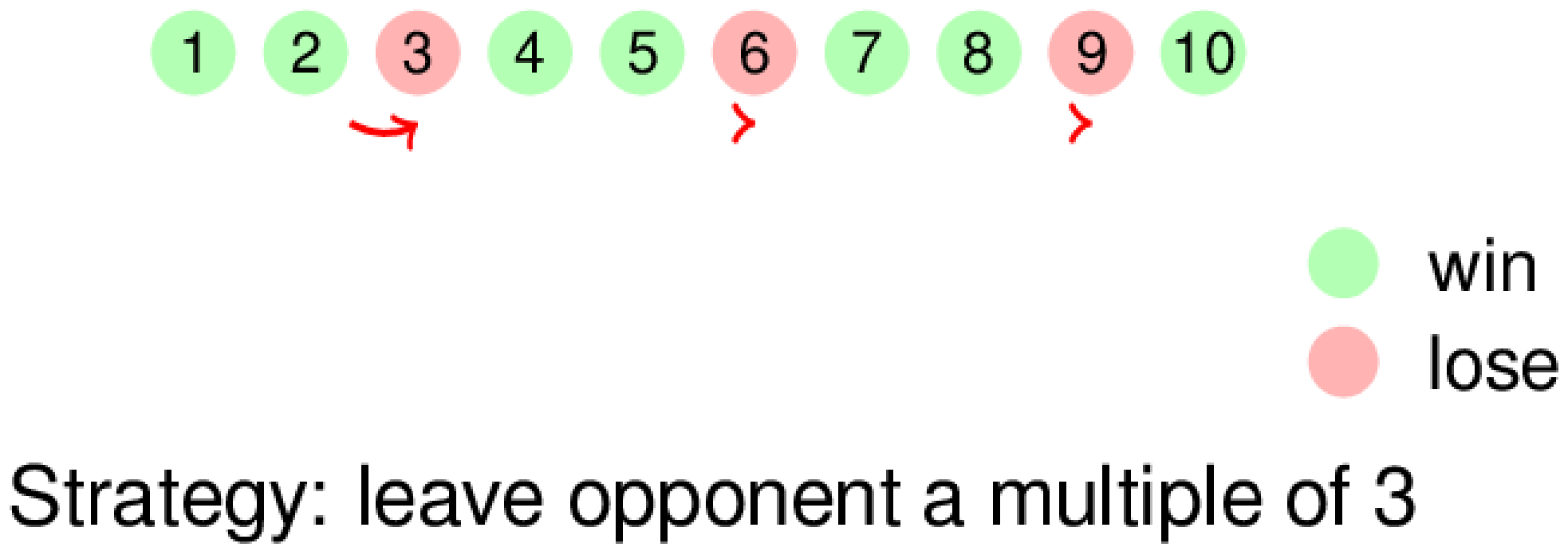

*Fig. 10. NIM*

### 7.1. NIM with one pile

A pile of **10 pebbles** (or candies). Two players take turns taking **1 or 2** pebbles. Whoever takes the **last** one wins.

The strategy was explained gradually, using **backward analysis** (backward induction):

**Step 1.** If **1** remains — you take it and win. □

**Step 2.** If **2** remain — you take 2 and win. □

**Step 3.** If **3** remain — you take 1 → opponent has 2 (they win). You take 2 → opponent has 1 (they win). **Losing position!**

**Step 4.** If **4** remain — take 1 → opponent has 3 (they lose!). □

**Step 5.** If **5** remain — take 2 → opponent has 3 (they lose!). □

**Step 6.** If **6** remain — whatever you take, opponent gets 4 or 5 and wins. **Losing position!**

Pattern: losing positions are **3, 6, 9, ...** — multiples of 3. Winning strategy: leave your opponent a number of pebbles that is a multiple of 3.

For a pile of 10: first move — take 1 (leave 9 = 3 × 3). Then after each opponent's move, complement to a multiple of 3.

### 7.2. Modification: 1, 2, or 3

When children mastered the "1 or 2" version, the rules changed: you can take **1, 2, or 3**. The method is the same, but the answer is different:

- Losing positions: **4, 8, 12, ...** — multiples of 4.
- Strategy: leave your opponent a multiple of 4.

Children had to **reconstruct** the new strategy themselves. Not receive a ready answer, but apply an already known method to new rules. General principle: if you can take from 1 to k pebbles, losing positions are multiples of (k+1). This wasn't stated to the children — they discovered the pattern through play.

### 7.3. NIM with two piles

Two piles of pebbles. Per turn, you can take as many as you want, but from only **one** pile. Whoever takes the last pebble wins.

Winning strategy: keep the piles equal in size. If the opponent takes from one — take the same amount from the other (mirror). Losing positions: (1,1), (2,2), (3,3), ... — when piles are equal.

### 7.4. Chodniczek

The name was invented by the children. A 10×10 board, the piece starts in the bottom-right corner. Per turn, you can move any number of squares, but either **left** or **up**. The goal is to reach the top-left corner. Whoever gets there wins.

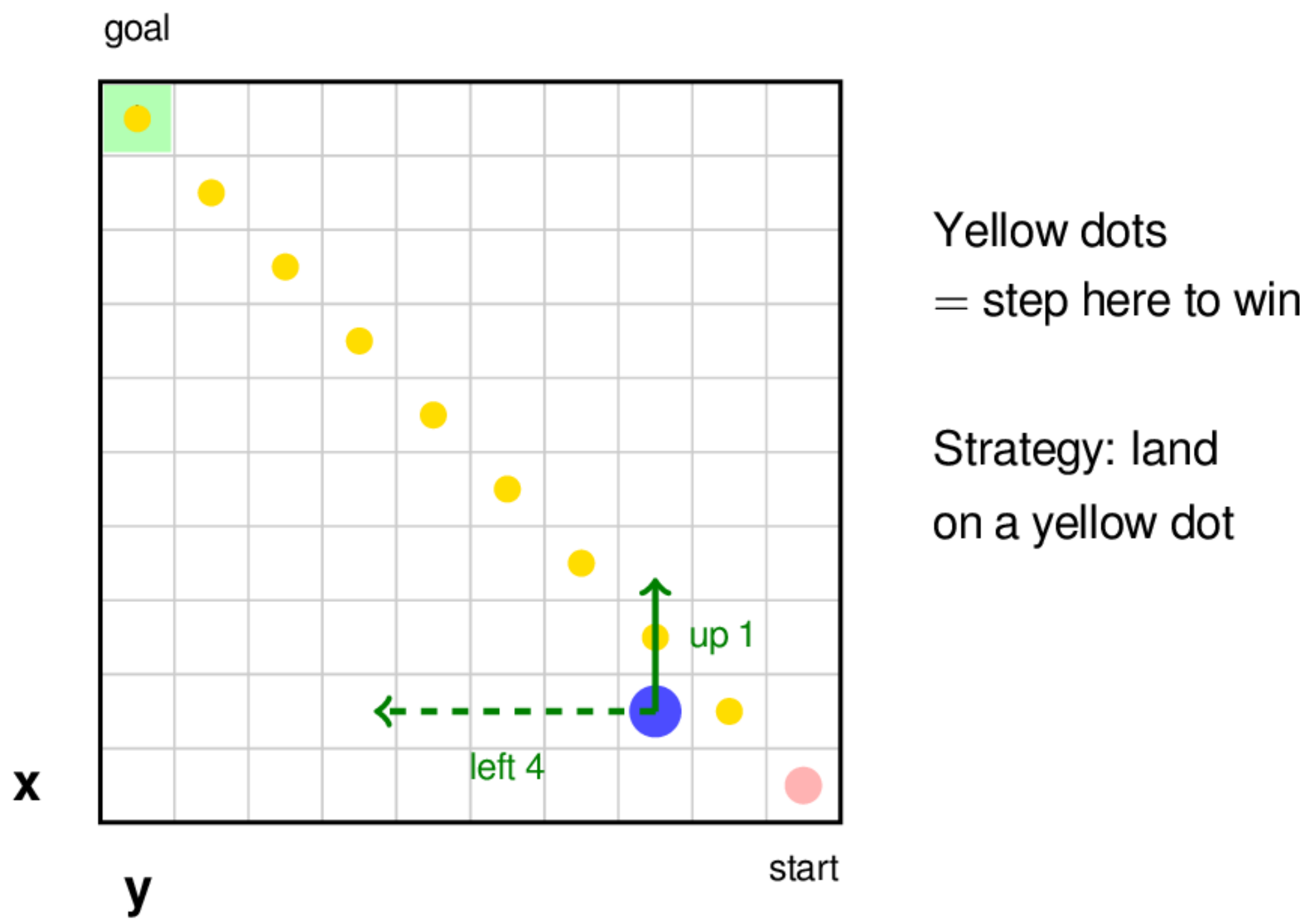

*Fig. 11. Chodniczek*

Winning strategy: step onto one of the yellow dots (the anti-diagonal, which in adult language is the set x+y=9). If you're on a yellow dot — any move takes you off, and the opponent can step back onto one. If you're not on a yellow dot — make a move onto one.

The board was usually 10×10, but sometimes we confused the children by changing the starting position — so the strategy couldn't be memorized but had to be understood.

### 7.5. Discovering isomorphism

These games were given at different times, and we didn't rush with explanations — as long as children were interested in simply competing. At some point we opened the Chodniczek table and asked: "*Which game we've played is this very similar to?*" The children didn't know.

Then we started writing side by side: when a child made a move in one game, we wrote the corresponding move in the other. Children **saw the correspondence themselves** — although we didn't explain it in words. They were able to match moves and transfer the winning strategy from one game to another.

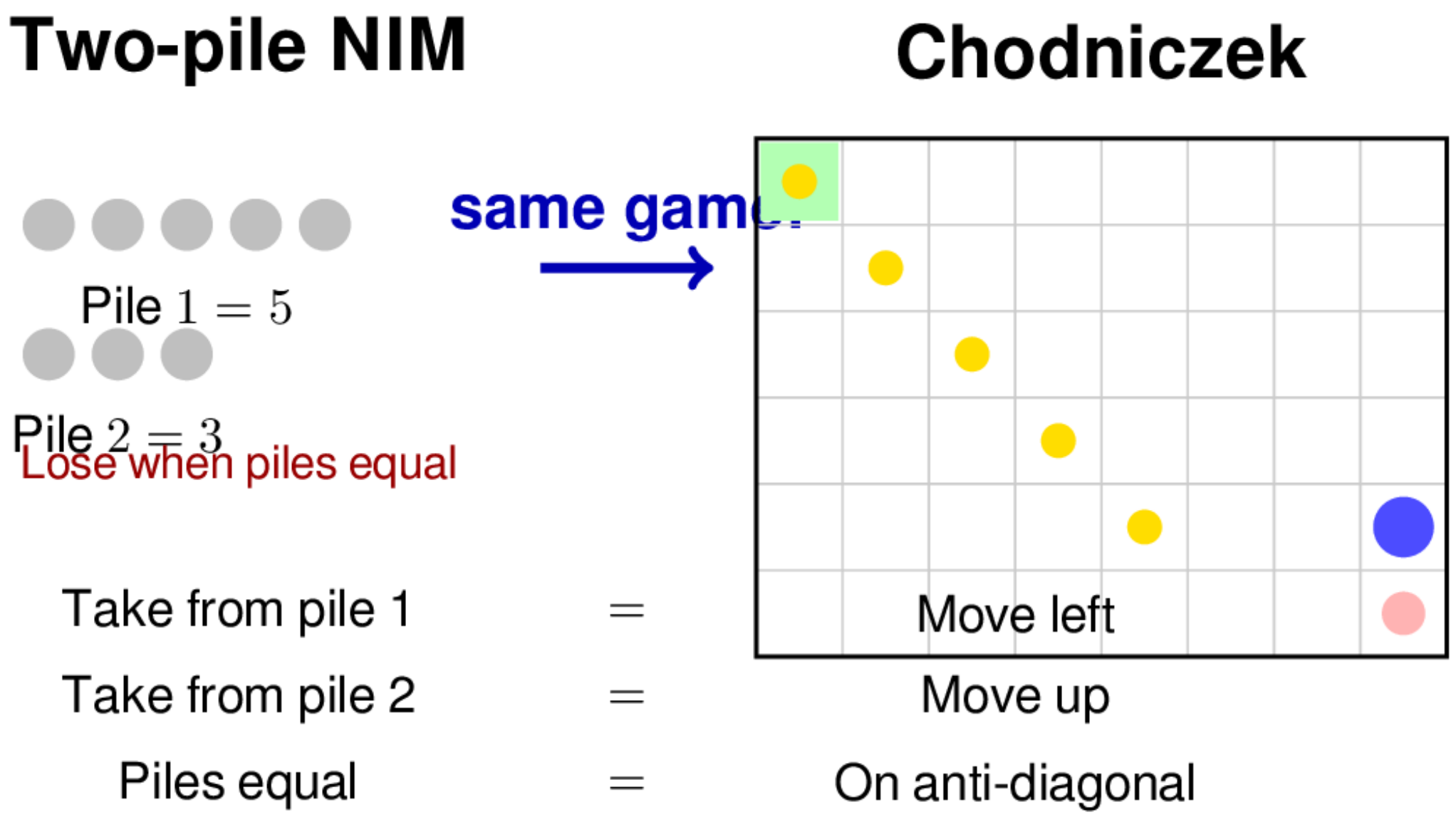

*Fig. 12. Isomorphism*

Position (x, y) on the Chodniczek board is the pair (pebbles in pile 1, pebbles in pile 2) in NIM. Moving left = taking from pile 1. Moving up = taking from pile 2. Anti-diagonal (equal piles) = loss.

This is one of the most valuable moments of the year: children saw that **two different games are the same**. In mathematics this is called isomorphism. Nine-year-old children don't know this word, but they experienced it.

## 8. Chess Miniatures

*Several sessions throughout the year.*

Another game that fit well into the competitive format — chess problems involving mate in a few moves. But with an important twist: one team plays for the strong side (must deliver mate in the given number of moves), while the other plays for the weak side (must survive and avoid mate). This equalizes chances: the team with the king doesn't need to know chess theory — just run away.

The figures below show two examples:

**Left:** "mate with two rooks in 4 moves": white king on d1, black rooks on b8 and h8. The rooks team must push the king to the edge of the board and deliver mate. The king team must make this as difficult as possible.

**Right:** "mate with queen and king in 5 moves": white queen on d3, white king on c1, black king on e5. The task is the same — deliver mate, but the queen alone is stronger than two rooks, so the plan is different: the king must be used to help.

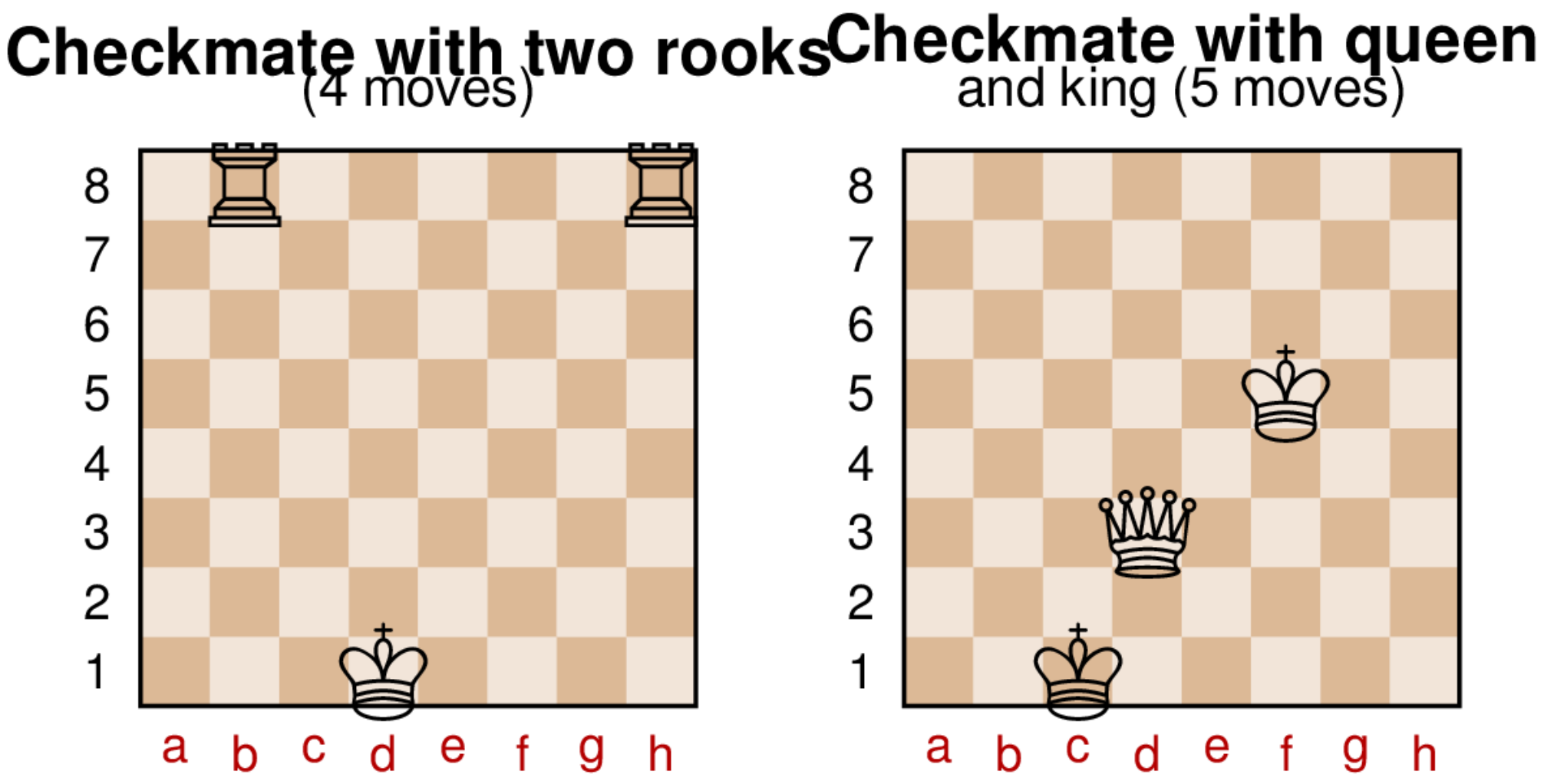

*Fig. 13. Chess miniatures*

The format is the same: a pair of students (one from each team) comes to the board, each makes a move for their side. The rest of the team can discuss, but the one at the board makes the move.

Chess problems develop a different type of thinking: you need to **plan several moves ahead**, accounting for the opponent's responses. This is the same backward analysis as in NIM, but in a much richer space — and with a visual, spatial component.

---

## 9. Arithmetic, the Equals Sign, and Scales

*A block of several sessions, arising from dissatisfaction with calculation speed in "Market."*

Playing "Market," we saw that children calculate slowly: multiplication and subtraction of two-digit numbers was difficult, 10 seconds for a decision often wasn't enough. We decided to devote separate sessions to arithmetic — but not as boring exercises, rather through a step-by-step subtraction method. In the process we discovered a deep problem with understanding the equals sign, which led to the idea of the scales lesson.

### 8.1. Step-by-step subtraction

Subtraction across ten: break the subtrahend into parts so the intermediate result lands on a round number.

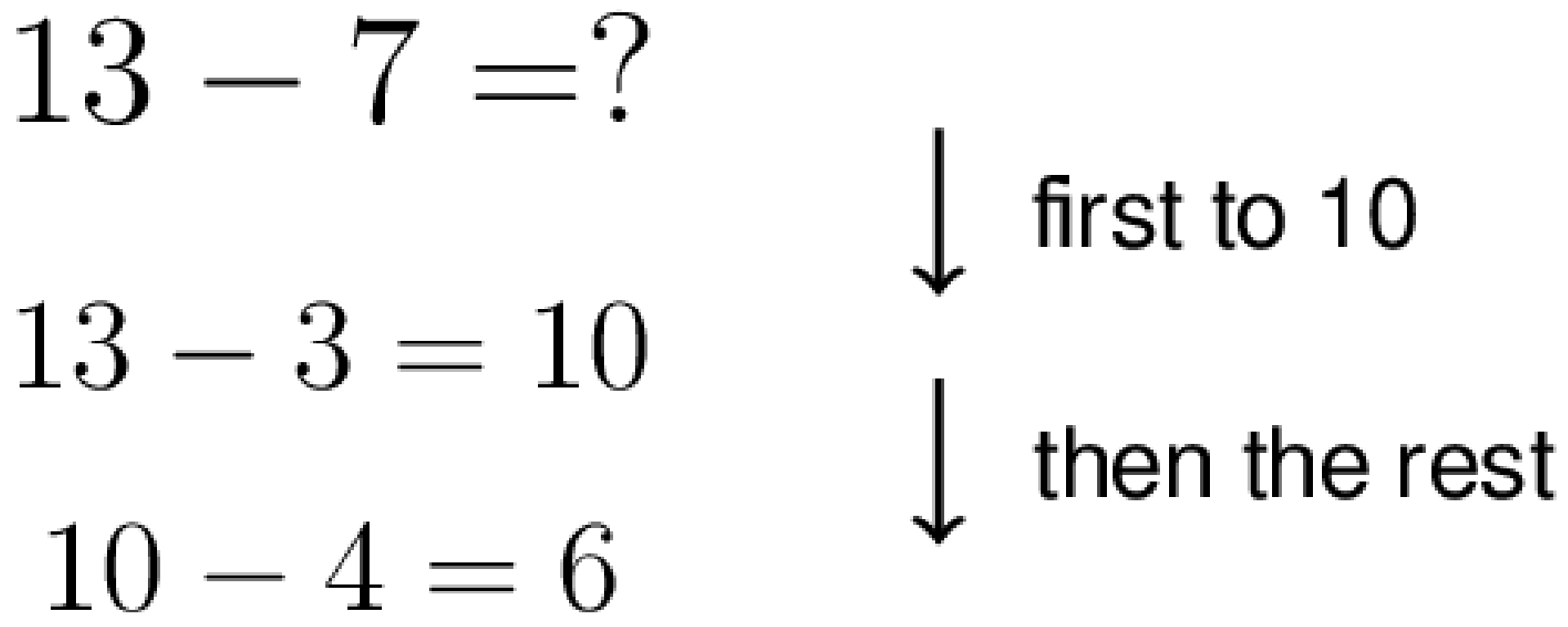

*Fig. 14. Subtraction across ten*

13 − 8 = 13 − 3 − 5 = 10 − 5 = 5

15 − 7 = 15 − 5 − 2 = 10 − 2 = 8

22 − 9 = 22 − 2 − 7 = 20 − 7 = 13

Children did exercises on paper — filling in the intermediate steps.

### 8.2. Discovery: children don't understand the equals sign

In the process it became clear that children **perceive the = sign as a mysterious symbol placed before the answer**. For them:

3 + 5 = 8

means: "three plus five, and here's the answer — eight." The = sign is an **arrow to the result**, not a statement about the equality of two expressions.

The notation 13 − 3 − 5 = 10 − 5 = 5 puzzled them: "Why are there several equals signs? Where's the answer?"

This is a widespread misconception described in the pedagogical literature. It hinders understanding of equations, algebra, and mathematical notation in general. Recognizing this problem led to the next part of the lesson — with scales.

### 8.3. Scales in the classroom

In the classroom, real balance scales with a set of weights were available. The lesson with them became one of the most memorable.

### 8.4. The merchant's problem

Children were told a story: a merchant traveled from village to village weighing goods. He carried 13 weights — from 1 to 13 kilograms. One day a clever person told him: "Why so many? Three weights are enough: **1, 3, and 9**. With them you can weigh any item from 1 to 13 kg."

### 8.5. The key idea: weights on both pans

To weigh **2 kg**, you need to place the item (2 kg) and the 1-weight on **one** pan, and the 3-weight on the **other**:

Left pan: 2 + 1 = Right pan: 3

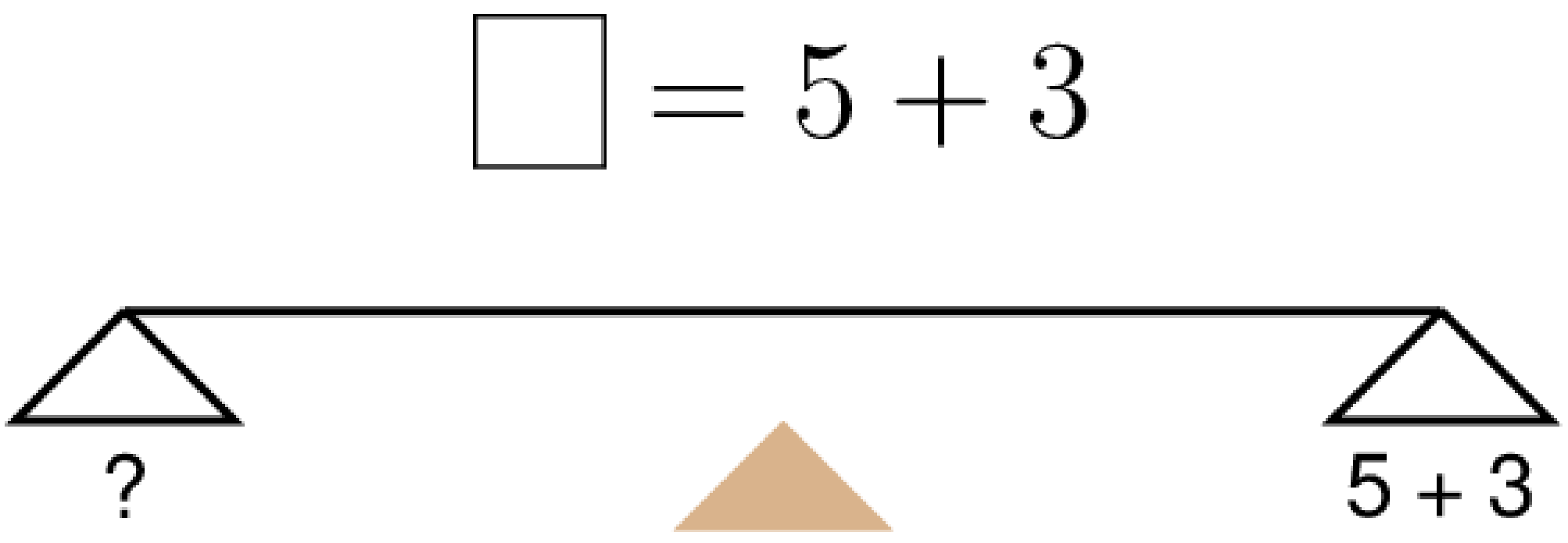

*Fig. 15. Scales*

This is impossible if weights go on only one pan. Children had to understand on their own that **weights can go on both pans** — otherwise even 2 kg can't be weighed.

## 8.6. All weights from 1 to 13

The authors showed the solution for 2 kg and wrote on the board: 2 + 1 = 3. Then for each subsequent number they asked: "Who's ready?" A child came to the board and wrote the equation.

| Weight (kg) | Left pan (item + weights) | Right pan (weights) | Notation |
|---|---|---|---|
| 1 | 1 | 1 | 1 = 1 |
| 2 | 2 + 1 | 3 | 2 + 1 = 3 |
| 3 | 3 | 3 | 3 = 3 |
| 4 | 4 | 3 + 1 | 4 = 3 + 1 |
| 5 | 5 + 1 + 3 | 9 | 5 + 1 + 3 = 9 |
| 6 | 6 + 3 | 9 | 6 + 3 = 9 |
| 7 | 7 + 3 | 9 + 1 | 7 + 3 = 9 + 1 |
| 8 | 8 + 1 | 9 | 8 + 1 = 9 |
| 9 | 9 | 9 | 9 = 9 |
| 10 | 10 | 9 + 1 | 10 = 9 + 1 |
| 11 | 11 + 1 | 9 + 3 | 11 + 1 = 9 + 3 |
| 12 | 12 | 9 + 3 | 12 = 9 + 3 |
| 13 | 13 | 9 + 3 + 1 | 13 = 9 + 3 + 1 |

## 8.7. Pedagogical effect

The = sign finally came to mean **equality** — the physical balance of the scales. The left pan weighs the same as the right. It's not "the sign before the answer" but a **fact you can see and touch**.

The notation 7 + 3 = 9 + 1 stopped being strange: the child sees on the scales that both pans are balanced. Both expressions are equal sides of an equality.

## 8.8. The mathematics behind the problem

The weights 1, 3, 9 are powers of three: $3^0$, $3^1$, $3^2$. Every weight from 1 to 13 can be represented as:

$w = a_0 \cdot 1 + a_1 \cdot 3 + a_2 \cdot 9$, where $a_i \in \{-1, 0, +1\}$

where +1 means "weight on the opposite pan," −1 means "weight on the same pan as the item," 0 means "weight not used." This is the **balanced ternary system**.

Children were not told this — but they were effectively working with it, selecting combinations.

---

## 10. Guess the Number

*Toward the end of the school year.*

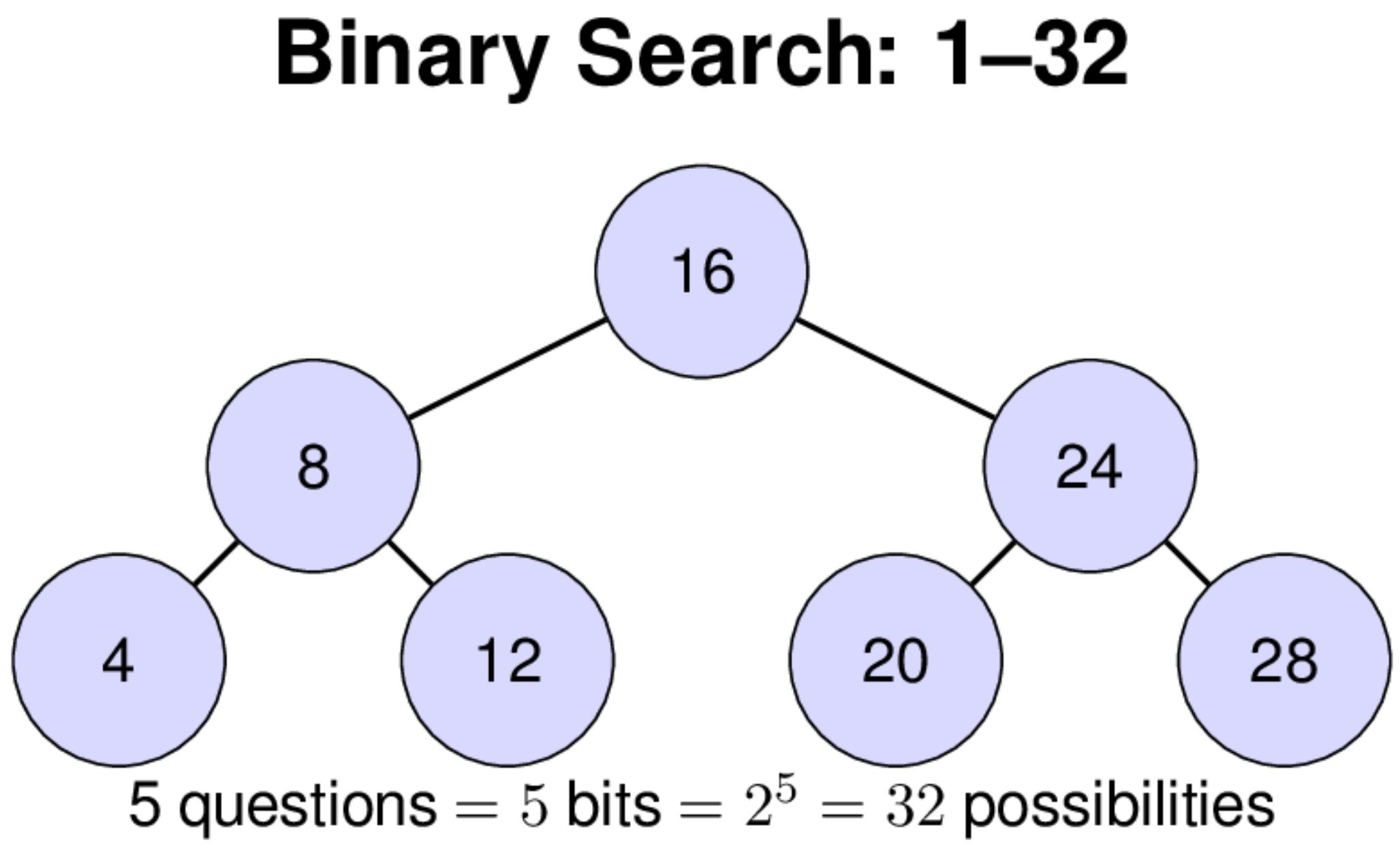

*Fig. 16. Binary search*

A player from one team thinks of a number from **1 to 32**. A player from the other team asks exactly **five** questions (answer only "yes" or "no"), then must name the number.

In practice, children asked questions of very different types. After each question, the teacher briefly wrote the question on the board and the list of remaining options — so the whole class could see how the set narrowed and which questions worked better.

Here's an example from an actual game (the number was 19):

| # | Child's question | Answer | Remaining on the board |
|---|---|---|---|
| 1 | Is the number even? | No | 1,3,5,7,9,11,13,15,17,19, 21,23,25,27,29,31 (16) |
| 2 | Does the number contain the digit 2? | No | 1,3,5,7,9,11,13,15,17,19, 31 (11) |
| 3 | Greater than 15? | Yes | 17,19,31 (3) |
| 4 | Greater than 20? | No | 17,19 (2) |
| 5 | Greater than 18? | Yes | 19 |

They guessed — but not optimally! The second question ("does it contain digit 2?") eliminated only 5 numbers out of 16 — unevenly. Children quickly noticed: **the best questions divide the remaining options in half**. The question "is it even?" splits 32 numbers evenly into 16+16, while "does it contain digit 2?" splits them unequally.

The board record made this visual: after a good question the list shrank by half, after a bad one it barely changed.

The optimal strategy is to halve the range:

**Example of optimal play.** The number is 19.

| # | Question | Answer | Range |
|---|---|---|---|
| 1 | Greater than 16? | Yes | 17–32 |
| 2 | Greater than 24? | No | 17–24 |
| 3 | Greater than 20? | No | 17–20 |
| 4 | Greater than 18? | Yes | 19–20 |
| 5 | Greater than 19? | No | 19 |

Five questions — five halvings: $32 = 2^5$. This is not an accidental choice of range.

Each yes/no question carries exactly **1 bit** of information. Five questions — 5 bits — can distinguish $2^5 = 32$ options. Children implicitly encounter **binary search**, **powers of two**, and **information as a resource** — a foolish question wastes the resource.

---

## 11. Summary

### What worked

• **Competition** — the main engine. Maintains the attention of the entire class, not just strong students.
• **Repeated play** of the same game with variations. Children gradually discover the strategy without receiving it "from above."
• **Minimal explanations.** Works in a language barrier situation, but useful in general: children learn through action, not through lectures.
• **Real objects** (scales, dice, broker stick) make abstract ideas tangible.
• **Pair work** mathematician + facilitator (Lakshtanov — content, problems, modifications; Tusiewicz — leads in Polish, Excel, materials; evening meetings).

### What we didn't get to

• **Traffic light races** — a game with speed and acceleration control, traffic lights with periodic switching.
• **Going to 2D** — races on a track (speed as a vector), space with gravity between planets.
• **Market with swaps.** An extension of the apple market introducing a new instrument — a deal for the future.

**Motivation — a story for children.** Before introducing the swap, we tell a story. Mr. Kowalski makes apple juice. He needs apples every month. If the price rises from 10 to 18, he'll go bankrupt — the juice will be too expensive. So he wants to agree in advance: "I'll buy apples from you next month at 13 złoty — regardless of the market price." The price of 13 is slightly above the current 10, but Mr. Kowalski is willing to pay this premium for peace of mind. And the farmer agrees: he gets 13 guaranteed, even if the price drops to 5. Both sides benefit — each has eliminated their risk.

**How it works in our game.** In our market, all teams are speculators — there is no separate "producer." Therefore swaps are sold and bought by the **facilitator** (Marek) — just as he sells and buys apples. At any point in the game, instead of a regular buy/sell, a team can ask: "We want a swap — buy 3 apples at 12 on turn 8." The facilitator writes the deal on the board. On turn 8 it executes automatically: the team pays $3 \times 12 = 36$ and receives 3 apples — regardless of the market price. All other information remains as before: on the projector screen children see the current price, each team's apples and money, the history of moves. The only addition is a line with the future deal on the board — a minimal change to the familiar format.

**Sample game with a swap:**

| Turn | Price | Event |
|---|---|---|
| 5 | 11 | Team A asks the facilitator for a swap: "Buy 3 apples at 12 on turn 8" |
| 5 | 11 | Facilitator writes on the board: "A ← 3 apples at 12, turn 8" |
| 6 | 14 | Price went up. Team A is pleased — will buy at 12 |
| 7 | 9 | Price dropped. Team A regrets — could have bought on the market at 9 |
| 8 | 13 | Swap executes: A pays 3×12 = 36 and gets 3 apples. Market price is 13 — A saved 3×(13−12) = 3 złoty |

If the price on turn 8 had been 10, team A would have **lost**: paid 36 instead of 30 on the market. A swap is not free insurance — it's a bet on the future.

**Why this is pedagogically sound:**

• Children already understand the market — the swap adds one new dimension (the future) without breaking the existing game.
• No need to explain the mechanism upfront. Just tell Mr. Kowalski's story, then say: "Now you can not only buy apples from the facilitator, but also arrange a future purchase. Want to try?" Children will learn through play.
• The question arises: **at what price** to enter the swap? The current price is 11, but a swap for turn 8 is a purchase in an unknown future. Should you pay more? Less? Children implicitly encounter the idea of **the price of risk**.
• The team must think **two moves ahead**: I have a swap on turn 8 — I need enough money by then. This adds complexity to planning, but within already familiar mechanics.

## Acknowledgements

The authors thank teachers Beata Lejman and Kamila Macieik Bania for their support and willingness to include the math circle in the class schedule. Special thanks to the students' parents for their warm feedback and support throughout the year.